\documentclass{article}

\usepackage{amsmath}
\usepackage{amssymb}
\usepackage{graphicx}
\def\be{\begin{equation}}
\def\ee{\end{equation}}
\def\bea{\begin{eqnarray}}
\def\eea{\end{eqnarray}}
\def\p{\partial}

\newcommand{\R}{{\mathbb R}}
\newcommand{\pp}{{\mathbb P}}

\title{Maximum compatibility estimates and shape reconstruction with
boundary curves and volumes of generalized projections
}

\author{Mikko Kaasalainen\\
Department of Mathematics\\
Tampere University of Technology\\
P.O. Box 553, 33101 Tampere\\ 
Finland}

\date{}
\begin{document}
\maketitle

\begin{abstract}
We show that the boundary curves (profiles) in $\R^2$ of the 
generalized projections of a body in $\R^3$
uniquely determine a large class of
shapes, and that sparse profile data, combined with projection
volume (brightness) data,
can be used to reconstruct the shape and the spin state of a body.
We also present an optimal strategy for the relative weighting of the data 
modes in the inverse problem, and derive the maximum compatibility 
estimate (MCE) that corresponds to the maximum likelihood or maximum
a posteriori estimates in the case of a single data mode. MCE is not
explicitly dependent on the noise levels, scale factors or numbers of data
points of the complementary data modes, and can be determined without the
mode weight parameters. We present a solution method well suitable
for adaptive optics images in particular, and discuss various choices of
regularization functions.
\end{abstract}


\section{Introduction}

Constructing the shape model of a three-dimensional object or a surface
is often based on images obtained at various viewing geometries. 
This is the standard case in human and robot vision as well as in 
cartography. When individual
points on the surface can be identified in different
images, one can directly solve for their
position vectors and use this as a top-down
basis for stereographic mapping or model construction. In many cases, 
however, the model construction reduces to an inverse problem 
rather than a cartographic one as the image information is based
purely on the projections of the target in the viewing directions.
When surface illumination and other effects are taken into account,
we talk about {\em generalized projections} \cite{genproj}
to distinguish these from
simple shadow-like projections or silhouettes. As discussed in \cite{genproj},
there are various types of such projections, ranging from the volume-like 
quantity of integrated brightness $L\in\R$
(generalization of the area of a shadow on a projection screen) to 
resolved images ${\cal I}\in\R^2\times\R$ (for one wavelength).

In this paper, we consider the case where images ${\cal I}(\omega,\omega_0)$
obtained at various viewing and illumination directions 
$\omega,\omega_0\in S^2$ are available, but the infomation in these
images is only contained in the boundary curves between the dark background 
or a shadow and the illuminated portion of the target surface. This situation
is typical for faraway objects in space for which low-resolution images
are available via large telescopes through adaptive optics (AO) \cite{keck} 
or other deconvolution and image processing techniques. Due to the
deconvolution process, the actual brightness levels of the pixels in these 
images tend to portray artificial and exaggerated features, so they 
are usually less reliable than profile contours, i.e.,
the locations of the light/dark boundary pixels \cite{carry1,carry2}. 

As the coverage of viewing geometries is seldom wide enough to enable a
full reconstruction of the model from images alone, we also
consider the possibility of augmenting the image dataset with a set of
measured brightnesses $L(\omega,\omega_0)$ of the target at
various observing geometries. Subsets of these, measured within some time 
intervals, are called lightcurves. As discussed in \cite{aa}-\cite{genproj} 
(see also references therein), a global, usually convex, model of the
target can be obtained from a large enough set of 
$L(\omega,\omega_0)$. The possibility to use images $\cal I$ both serves
to reconstruct more details in the model and to use a combined dataset
for successful modelling when neither the available $L$ nor $\cal I$ 
are sufficient alone.

The ratio of profile contour pixels to all pixels of the target is
approximately $4/D$, where $D$ is a typical diameter of the target in
pixels. For the largest few asteroid AO targets, $D$ is around 30. When $D$
is less than about 10, the location accuracy of 
the border pixels is not necessarily very
much better than the brightness accuracy of the pixels,
but on the other hand there is no particular loss in the number of
information points when only borders are used. For higher $D$, border points 
are a smaller subset of available pixels, but now their location accuracy is 
far better than the brightness accuracy of all pixels \cite{carry1,carry2}. 

The profile contours of the AO images are obtained as a solution of
a separate inverse (or imaging) problem, where an approximation of the
atmospheric point-spread function (PSF) is 
first used to deconvolve the raw image
(with, e.g., connectedness of the processed image as a constraint),
and the contours can be separately modelled with wavelet techniques 
\cite{carry1,carry2,keck}. This independence from the actual model of
the target is advantageous in the sense that the assumptions and
inevitable deficiencies of the model (particularly in the adopted
light-scattering model on the surface \cite{iau}) do not affect the
outcome of the AO image processing. 
On the other hand, from the methodological point of view, all information 
should usually be employed simultaneously when solving an inverse problem, 
so another approach would be to use the raw AO image data and the 
approximated PSF directly in model construction without separate image 
deconvolution. However, in practice it appears that the profile curve
extraction procedure, in particular, retains valuable independent 
information \cite{carry1,carry2}, 
and the model deficiencies affect the fit deviation between
the predicted and observed model profile 
curves less than they affect that between the full model and AO images.
What is more, below we show that the profile
curves convey almost as much information on the shape as the full images,
so we can conclude that the two-step inversion of AO (and brightness) 
data is well justified.

The paper is organized in sections in the following manner. In Section 2
we study the information content of profile contours and show uniqueness
results for the inverse problem of reconstructing shapes from these.
Section 3 deals with the posing of the inverse problem and the
choice of regularization functions. In Section 4, we discuss the
weighting and maximum compatibility estimates for inverse 
problems with multiple data modes,
and examples of the use of real brightness data $L$ and
images $\cal I$ are presented in Section 5. Section 6 sums up.

\section{Generalized profiles: uniqueness results}

In this section, we define the concept of generalized profiles as 
boundary curves corresponding to generalized projections, and show that
a large class of shapes is reconstructable from these. While many practical
procedures for shape-from-profiles (also known as volume carving) and 
shape-from-shading are well known in, e.g., computer vision 
(see \cite{shadcar} and references therein) and 
cartography (clinometry), some of their
geometric characteristics and the properties of the
corresponding shape classes and inverse problems 
discussed here have not been 
previously stated in the mathematical literature, to the
best of our knowledge.

Let us first look at
classical profiles defined by one direction $\omega\in S^2$, i.e.,
$\omega_0=\omega$.
The projection ${\cal P}(\omega,{\cal B})\in\R^2$ of a compact set 
${\cal B}\in\R^3$ (a set of points on closed surfaces) 
in the direction $\omega\in S^2$
maps $x\in{\cal B}\rightarrow \varkappa\in\cal P$:
\be
\varkappa=\left(\begin{array}{ccc} 0&1&0\\
0&0&1\end{array}\right){\sf R}_y(\vartheta-\frac{\pi}{2}){\sf R}_z
(\psi)\, x,
\ee
where $0<\vartheta<\pi,0\le\psi<2\pi$ are the polar coordinate
angles defining 
$$
\omega=(\sin\vartheta\cos\psi,\sin\vartheta\sin\psi,
\cos\vartheta),
$$ 
and ${\sf R}_i(\alpha)$ is 
the rotation matrix corresponding to the rotation 
of the coordinate frame through angle $\alpha$
in the positive direction about the $i$-axis. 
Thus, e.g., ${\sf R}_z(\phi)$ is
\be
{\sf R}_z(\phi)=\left(\begin{array}{rrr}\cos\phi&\sin\phi&0\\
-\sin\phi&\cos\phi&0\\
0&0&1\end{array}\right).\label{rz}
\ee
For $\vartheta=0$, $\varkappa_i=x_i,i=1,2$, and for
$\vartheta=\pi$, $\varkappa_1=-x_1,\varkappa_2=x_2$.
Only a half of $S^2$ is needed for defining $\omega$ as
$$
\varkappa_i(\omega)=(-1)^i\varkappa_i(-\omega).
$$

For some definitions below, we need to give the 
three-dimensional position
$x_\varkappa(\omega)\in\R^3$ of the planar points $\varkappa$:
\be
x_\varkappa(\omega)={\sf R}_z(-\psi){\sf R}_y(\frac{\pi}{2}-\vartheta)\left(
\begin{array}{c}0\\ \varkappa_1\\ \varkappa_2\end{array}\right).
\ee

{\bf Definition 1.}
The {\em profile} ${\cal P}_\p (\omega,{\cal B})$ of $\cal B$ in the
direction $\omega$ is the boundary of its projection 
${\cal P}(\omega,{\cal B})$:
$$
{\cal P}_\p(\omega,{\cal B})=\p{\cal P}(\omega,{\cal B})
$$
(the specific notation ${\cal P}_\p$ is used for emphasis).
More specifically, ${\cal P}_\p (\omega,{\cal B})$
is the mapping $x\rightarrow\varkappa$ from those points $x\in\cal B$
for which, for all lines $x+s\omega$ parametrized by $s$,
\be
{\cal P}_\p (\omega,{\cal B})
=\left\{ \varkappa\Big\vert\forall s:\exists\theta,\exists\rho>0:
{\mathfrak P}(x+s\omega;\theta,\varrho)\notin{\cal B}, 
0<\varrho<\rho\right\},
\ee
where ${\mathfrak P}({\cal L};\theta,\varrho)$ 
denotes a parallel transport, perpendicular to the line $\cal L$, of $\cal L$
by the amount $\varrho$ in the direction $\theta\in S^1$ in some system
around $\cal L$.

{\bf Definition 2.}
The {\em cylinder continuation} (CC) ${\cal C}(\omega,{\cal S})$
of a set $\cal S$ of points $x\in\R^3$ is the set of
points in $\R^3$ given by
\be
{\cal C}(\omega,{\cal S})=\Big\{x+s\omega\Big\vert 
x\in{\cal S};-\infty<s<\infty\Big\}.
\ee
 
{\bf Definition 3.} The {\em profile hull} ${\cal H}\in\R^3$ is the
bounding surface of the 
set of points formed by the intersection of the CCs of the projections 
$\cal P$ in
some directions $\omega_i$, $i=1,\ldots,N$, corresponding
to measured profiles ${\cal P}_\p(\omega_i)\in\R^2$:
\be
{\cal H}(\{{\cal P}_\p (\omega_i)\vert_{i=1}^N\})
=\p\bigcap_i {\cal C}(\omega_i,{\cal S}_i),\quad
{\cal S}_i=\Big\{ x_\varkappa(\omega_i)\Big\vert\varkappa
\in{\cal P}(\omega_i)\Big\}.
\ee

{\bf Remark.} For practical purposes and clarity, we assume $\cal H$
to be constructed such that it has a 
closed surface as a boundary, though this is not strictly
necessary in its
definition (we could define $\cal H$ as set of points rather than its 
bounding surface). Thus $\cal B$ and $\cal H$ are in the same object class.
Similarly, we use the concepts of projection and profile (or body
and surface) somewhat
interchangeably when the meaning is obvious.

Many bodies can be reconstructed to arbitrary accuracy
with profile hulls.
A convex body is already 
determined by its $\cal H$ constructed with a full coverage
of the directions $\omega$ confined to any plane in $\R^3$.
In real profile measurements, however,
the position of the profile in the $\varkappa$-plane
is usually arbitrary, i.e., the profile is determined up to a translation
$\varkappa_0$ of the profile plane origin. Then data restricted to 
planar directions are not necessarily sufficient even for convex bodies:
the profile offsets $\varkappa_{0,i}$ of a convex body are not always
uniquely defined by the profiles ${\cal P}_\p (\omega_i)$ via the profile 
hull $\cal H$ when $\omega_i$ are confined to a plane. This is simple to
illustrate by
considering curves in $\R^2$ and their projections in $\R$ at directions
in $S^1$.

If the curvature function $C:S^1\rightarrow\R$ 
of a closed convex curve in $\R^2$ is given
by a real-valued non-negative Fourier series
$$
C(\varphi)=\Re\sum_n c_n e^{in\varphi}\ge 0,\quad n\ge 0,
$$
where $\varphi$ denotes the direction of the outward normal of the curve,
the projected width $w(\varphi)$ of the contour in that direction is
$$
w(\varphi)=\int_{-\pi/2}^{\pi/2}\,C(\psi+\varphi)\cos\psi\,d\psi
=\Re\sum_n c_n e^{in\varphi}\,I_n,
$$
where
$$
I_n=\int_{-\pi/2}^{\pi/2}\cos n\psi\cos\psi\,d\psi
$$
and for $n\ne1$
$$
I_n=\left\{ \begin{array}{rl}
0&n=3,5,7,\ldots\\
2/(n^2-1)&n=2,6,10,\ldots\\
2/(1-n^2)&n=0,4,8,\ldots
\end{array}\right.
$$
and $c_1\equiv 0$ since $I_1\ne0$ and we must have 
$w(\varphi)=w(\varphi+\pi)$. Thus $w(\varphi)$ carries no information
on the odd-valued $n$-coefficients of the curvature function $C(\varphi)$
that uniquely defines the shape of the curve (cf.\ \cite{aa} for convex
surfaces in $\R^3$).
Thus, for cylindrical convex surfaces in $\R^3$, the observed profiles
in the symmetry plane can be made to correspond to any odd parts of
$C(\varphi)$ with suitable chosen offsets $\varkappa_0$. A typical case
is that of shapes mimicking 
a circular cylinder with constant $C$ in the symmetry plane:
for example, if we change $\varkappa_0=0$ to
\be
\varkappa_0^{(1)}(\varphi)=R(-1)^{\varphi\,{\rm div}\frac{\pi}{3}}
[\frac{2}{\sqrt 3}\cos(\varphi\,{\rm mod}\frac{\pi}{3}-\frac{\pi}{6})-1],
\ee
the reconstructed shape is a cylindrical Reuleaux triangle. This degeneracy
occurs since for projections $\R^2\rightarrow\R$ the volume of a profile 
is equivalent to the profile itself up to an offset. An additional
profile from a direction perpendicular to the plane resolves the degeneracy.

Thus, in general, we need data at full $\omega\in S^2$ for a unique
reconstruction of a body when the profile offsets are not known.
The principle in the reconstruction via the profile hull 
$\cal H$ is the requirement that
$\cal H$ must be consistent with the observed profiles,
i.e., the profiles of the constructed $\cal H$  must be identical
to the observed ones:
$$
{\cal P}_\p[\omega_i,{\cal H}(\{{\cal P}_\p (\omega_j)\vert_{j=1}^N\})]
={\cal P}_\p (\omega_i)
$$ 
(otherwise the
volume of the intersection defining $\cal H$ is not maximal).

Let us denote by the
{\em complete profile hull} ${\cal H}_C$
the profile hull for which $\omega$ covers all of $S^2$. The
complete profile hull ${\cal H}_C({\cal B})$ of a body $\cal B$
is the envelope of those of its 
tangents that do not intersect $\cal B$ anywhere. Then we can
define a class of surfaces that includes convex ones but extends far
into nonconvex surfaces:

{\bf Definition 4.}
{\em Tangent-covered bodies} (TCBs) are bodies that are their own
complete profile hulls: ${\cal B}={\cal H}_C({\cal B})$. Thus,
each surface point $x\in{\cal B}$ of a TCB is mapped at least to one
${\cal P}_\p(\omega)$. 

TCBs include a large variety of nonconvex surfaces or sets
of them: for example,
a body consisting of two separate spheres is a TCB. While convex bodies
$\cal C$ are reconstructable from the volumes of their generalized 
projections \cite{aa,genproj} and are
defined by having no tangents intersecting the body, TCBs $\cal T$
are reconstructable
from profiles and are defined by there being at least one tangent at
each surface point not intersecting the body elsewhere.

Let us now generalize the concept of profile in the same way as
projections. This leads to a shape class $\cal G$,
larger than TCBs $\cal T$, that can be reconstructed from generalized
profiles:
$$
{\cal C}\subset{\cal T}\subset{\cal G}.
$$

When we consider the directions $(\omega,\omega_0)$ in $S^2\times S^2$,
the region both visible and illuminated on $\cal B$ is given by 
\cite{genproj, kaast}
\be
{\cal A}_+(\omega,\omega_0;{\cal B})={\cal A}_+(\omega;{\cal B})\cap
{\cal A}_+(\omega_0;{\cal B}),
\ee
where
\be
{\cal A}_+(\omega;{\cal B})=\left\{x\in{\cal B}\Big\vert
\langle\nu(x),\omega\rangle\ge 0; \forall s>0: x+s\omega\notin{\cal B}
\right\},\label{aplus}
\ee
where $\nu(x)$ is the unit surface normal at $x$. The projection $\cal P$
of the boundary $\p{\cal A}_+$ is now the generalized profile:

{\bf Definition 5.} The {\em generalized profile} of the body $\cal B$
in the direction $\omega$ and at illumination direction $\omega_0$ is
\be
\p{\cal P}[\omega,{\cal A}_+(\omega,\omega_0;{\cal B})]=
{\cal P}[\omega,\p{\cal A}_+(\omega,\omega_0;{\cal B})].
\ee

The shape class $\cal G$ is not as straightforward to define as $\cal C$ and
$\cal T$. We can, however, prove configurations allowing unique shape
determination that illustrate its extension from 
$\cal T$.

{\bf Theorem 1.} {\em Assume that we know some parts $\cal K$
of a body $\cal B$ from profile measurements. There exist 
configurations in which unknown parts $\cal U$ of $\cal B$ not determinable 
from profiles can be uniquely
determined from the generalized profiles of $\cal B$ by using
the shadow boundaries of $\cal K$ on $\cal U$.}

{\em Proof.}
Assume that ${\cal H}_C({\cal B})$ is defined, and that it contains a 
planar section $\cal K$, and that $\cal B$ contains in this region
an unknown concavity ${\cal U}$ (corresponding to 
${\cal H}_C\setminus{\cal B}$). Also, assume that all points of $\cal U$
are seen from the viewing direction $\omega\bot\cal K$, and that the
illumination direction $\omega_0$ lies in a plane $\bot\cal K$, with
$\theta=\measuredangle(\omega_0,{\cal K})$. Then the planar edge curve
of the concavity $\cal U$
can be determined when $\theta\rightarrow 0$:
$$
\p{\cal U}=\lim_{\theta\rightarrow 0}\p S(\theta):=\p S_0,
$$
where $\p\cal S$ denotes the projection of the shadow boundary in the 
direction $\omega$ on the plane $\cal K$.
At various $0<\theta\le\pi/2$, we can measure the shadow boundary projections 
$\p{\cal S}(\theta)$, and thus extract the projection 
$\p\tilde S(\theta)$ of the shadow boundary inside $\cal U$:
$$
\p\tilde S(\theta)=\p S(\theta)\setminus\p\hat S(\theta),\quad 
\p\hat S(\theta):=\p S(\theta)\cap\p S_0.
$$
Then the envelope in $\R^3$ of the
intersection curves of the cylinder continuations of
$\p\tilde S(\theta)$ in the directions of $\omega$
and $\omega_0$
$$
{\cal C}[\omega,\p\tilde S(\theta)]\cap {\cal C}[\omega_0,\p\hat S(\theta)]
$$
uniquely constructs the surface of the concavity $\cal U$ (when $\cal U$
is suitably regular). $\square$

{\bf Theorem 2.} {\em There exist configurations
in which unknown parts $\cal U$ of $\cal B$
can be uniquely determined by using their shadow boundaries on the known part
$\cal K$.}

{\em Proof.}
Let $\cal B$ be a combination of a TCB $\cal E$ and any surface $\cal D$ 
(${\cal E}\cap{\cal D}=\emptyset$) that
can be determined using profiles in the directions $\omega$
for which the profile intersection
$$
{\cal Q}(\omega):={\cal P}(\omega,{\cal D})\cap{\cal P}(\omega,{\cal E})
$$
vanishes, $\cal Q=\emptyset$ (the whole of $\cal D$ is in the known
part $\cal K$). The unknown parts are assumed to be
on $\cal E$ (they cannot be determined using the above $\omega$).
Now the parts of profiles of $\cal E$ that merge with 
${\cal P}(\omega,{\cal D})$ at some $\omega$, i.e., 
${\cal Q}(\omega)\neq\emptyset$, are represented as shadows on $\cal D$
that we assume we can see from some directions $\omega'$.
The full or partial profiles of $\p{\cal P}(\omega,{\cal E})$ for which
$\p{\cal P}(\omega,{\cal E})\cap {\cal Q}=\emptyset$ can be determined
as usual, and, with a known $\cal D$,
the remaining parts $\p{\cal P}(\omega,{\cal E})\cap {\cal Q}\neq\emptyset$
can be determined from the shadows on $\cal D$. The intersection
$$
\p{\cal W}={\cal C}[\omega',\p S_p(\omega')]\cap{\cal D},
$$
where ${\cal C}[\omega',\p S_p(\omega')]$ denotes the cylinder continuation
corresponding to the observed projection of the shadow boundary of $\cal E$
on $\cal D$ in the direction $\omega'$, can be used to determine the
projection ${\cal P}(\p{\cal W},\omega)$, which completes the missing parts
of the needed profiles. Now we have 
constructed the set of full profiles of $\cal E$ at all $\omega\in S^2$, so
$\cal E$ can be constructed as it is a TCB. $\square$

Continuing in a similar manner, we can construct more complex 
variations of the two cases above to explore the shape class $\cal G$.
In practice, directions $(\omega,\omega_0)$ seldom cover $S^2\times S^2$
extensively or densely, so the shape is reconstructed within some
resolution (discretization degree of the model) and a priori assumptions,
as discussed below.

\section{Inverse problem}

Let us now consider the inverse problem of determining the shape and spin state
of a body $\cal B$ from some measured generalized profiles 
$\p{\cal P}[\omega_i,{\cal A}_+(\omega_i,\omega_{0i};{\cal B})]$
and the volumes $L(\omega_{0i},\omega_i)$ of generalized projections.
We present a method that is suitable for typical ground-based astronomical
data, i.e., the profiles are only obtained at restricted geometries and
their resolution level is not high. When a dense coverage of geometries
and high resolution are available (e.g., space probe missions), direct
methods of computer vision and cartography are usually applicable. 

Our goal is to construct a total goodness-of-fit measure $\chi^2_{\rm tot}$ 
\be
\chi^2_{\rm tot}=\chi^2_L+\lambda_\p\chi^2_\p+\lambda_R g(P),
\ee
where $L$ denotes lightcurves, $\p$ generalized profiles, and $R$ 
regularizing functions $g(P)$, where $P\in\R^p$ is the vector of
model parameters. Regularization is discussed at the end of this 
section, and the determination of the weights $\lambda$ in section 4.
We note here that the additional $g(P)$ make
$\chi^2_{\rm tot}$ pseudo-$\chi^2$ as it no longer describes an underlying
(assumed) Gaussian probability distribution
(though $g(P)$ may be $\chi^2$-like in their functional form). When
probability densities such as a posteriori distributions are constructed
from $\chi^2_{\rm tot}$, one can
assume $\chi^2$-distributions (of the form $e^{-c\chi^2}$) only for the
data components, and other suitable (prior) distributions for
the regularization components \cite{kaipio} such that the maximum
of the a posteriori distribution occurs at $\arg\min\chi^2_{\rm tot}(P)$.

Throughout this paper, we do not include the conventional 
$1/\delta^2$-factor in $\chi^2$-forms, where
$\delta$ is the expected (Gaussian) error variance (noise level), since
$\delta$ is seldom known exactly, and it does not affect the determination
of our point estimates which is the goal of this paper. 
Suitable parameters for Gaussian or other
distribution widths can be inserted separately whenever 
we want to construct a distribution.

The volumes of generalized projections are also called total or
disk-integrated brightnesses \cite{aa,kaast,genproj}:
\be
L(\omega_0,\omega)=\int_{{\cal A}_+} 
R(x;\omega_0,\omega)\langle\omega,\nu(x)\rangle\,  d\sigma(x)
\equiv \int_{{\cal P}(\omega,{\cal A}_+)} 
R[P^{-1}(\omega,{\cal A}_+,\varkappa);\omega_0,\omega] 
d^2\varkappa,
\ee
where $\nu(x)$ and $d\sigma(x)$ are, respectively, the
outward surface normal and surface patch of $\cal B$,
$R(x;\omega_0,\omega)$ describes the intensity of 
scattered light at the point $x$ on the surface,
$P^{-1}(\omega,{\cal A}_+,\varkappa)$ is the point in ${\cal A}_+$
corresponding to the projection point $\varkappa$,
and $d^2\varkappa$ is the surface patch of the projection $\cal P$.
In its basic form,
\be
\chi^2_L=\sum_{i} [L^{\rm (obs)}(\omega_{0i},\omega_i)
-L^{\rm (mod)}(\omega_{0i},\omega_i)]^2\label{lcchisq}
\ee
(assuming a constant noise level; 
see \cite{kaast} for modifications and variations of this). $L$-data
on $S^2\times S^2$ uniquely determine a convex body and the
solution is stable \cite{aa,genproj}, but $L$-data do not carry information
on nonconvexities in most realistically available $S^2\times S^2$ 
geometries in practice \cite{durech}. Such information must be provided
by AO or other techniques.

For many typical AO targets, the generalized profiles are starlike
due to the proximity of $\omega$ and $\omega_0$ and some regularity
of the target shape at global scale \cite{carry1,carry2,keck}. Then 
we can write $\chi^2_\p$ by considering, for each profile $i$,
their observed and modelled
maximal radii (from some point within the profile) 
on the projection plane $\varkappa=(\xi,\eta)\in\R^2$ at direction angles 
$\alpha_{ij}$ (starting from a chosen
coordinate direction for positive $\xi$, $\eta=0$):
\be
\chi^2_\p=\sum_{ij} [r_{\rm max}^{\rm (obs)}(\alpha_{ij})
-r_{\rm max}^{\rm (mod)}(\alpha_{ij})]^2.\label{aochisq}
\ee
As the accuracy of $r_{\rm max}^{\rm (obs)}$ is proportional to the size
of the image, the sum (\ref{aochisq}) automatically takes this weighting
into account (of course, direct weighting due to varying
noise levels can be used as well).

We now represent the body $\cal B$ as a polytope \cite{kaast}. Let
two vertices $a$ and $b$ of a facet have projection points $(\xi_a,\eta_a)$,
$(\xi_b,\eta_b)$, and $(\xi_0,\eta_0)$ be the point on the projection plane
from which the radii and $\alpha$ are measured (this defines the profile
offset that must be solved for in the inverse problem). With
\be
\begin{array}{ll}
A=-\sin\alpha, &B=\cos\alpha,\\
C=\eta_a-\eta_b, &D=\xi_b-\xi_a,\\
E=A\xi_0+B\eta_0, &F=\xi_b\eta_a-\xi_a\eta_b,
\end{array}
\ee
the intersection point of the radius line and the projection of the
facet edge $ab$  is at
\be
\xi=\frac{DE-BF}{AD-BC},\qquad \eta=\frac{AF-CE}{AD-BC},\label{intersect}
\ee
and, to be in the correct direction of $\alpha$ and between the points
$a$ and $b$, the intersection point must satisfy
\bea
(\xi-\xi_a)(\xi_b-\xi)\ge0,\quad &(\xi-\xi_0)\cos\alpha\ge 0,\\\nonumber
(\eta-\eta_a)(\eta_b-\eta)\ge 0\quad &(\eta-\eta_0)\sin\alpha\ge 0.
\eea
If $AD-BC=0$, the line in $\alpha$-direction
is parallel to the line $ab$, so there is no intersection unless
the lines coincide, i.e., either of the numerators in (\ref{intersect})
vanishes.

The model $r_{\rm max}(\alpha)$ 
can now be determined by going through all eligible
facet edges and their intersection points $\varkappa_{ab}(\alpha)$:
\be
r_{\rm max}^{\rm (mod)}(\alpha)=\max\Big\{\Vert
\varkappa_{ab}(\alpha)-\varkappa_0\Vert
\Big\vert a,b\in {\cal V}_+\Big\},
\ee
where ${\cal V}_+$ is the set of vertices of the facets ${\cal A}_+$ that are
both visible and illuminated. The set ${\cal A}_+$ of (\ref{aplus})
is determined by ray-tracing
\cite{kaas01}. It is an approximation (correct to the order of the 
average facet area) of the actual visible and illuminated
region, i.e., each facet either is or is not in ${\cal A}_+$ (judging by its
centroid): projection lines of
obstructing facets inside a facet are neglected when the facets are small 
enough.

The principle of using outer contours applies to AO data that do not generally
show non-starlike or multiple contours (due to crater shadows)
as the solar phase angles $\arccos\langle \omega,\omega_0\rangle$
are low and the resolution/deconvolution 
accuracy is not high. At high phase angles, even starlike bodies form
non-starlike contours, and contours inside the outer contour appear in
high-resolution images from probe flybys.

The outer contour $\p\cal O$
can be automatically derived for non-starlike shape models
as well; such models can be constructed by, e.g., joining starlike
submodels together, using a cylindrical coordinate frame \cite{kaast}, 
or by determining the coordinates of a set of points with which a
suitable surface (a new tessellation for each iteration) 
is defined via, e.g., mesh-free methods
such as weighted/moving least squares \cite{levin}. For clarity,
let us first assume that no other generalized profile contours exist
outside $\p\cal O$. Denoting the edges of the facets
of ${\cal A}_+$ by ${\cal E}_+$,
$\p\cal O$ is constructed by the following algorithm:

{\bf 1.} Construct the set ${\cal F}_0\subset {\cal E}_+$
of the edges of ${\cal E}_+$ 
that are shared by a facet in ${\cal A}_+$ and by
a facet not in ${\cal A}_+$
but for which $\langle \nu,\omega \rangle\ge 0$.

{\bf 2.} Construct the set
${\cal F}\subset {\cal E}_+$ of the edges shared by a facet of ${\cal A}_+$
and a facet
for which $\langle \nu,\omega \rangle<0$.

{\bf 3.} Construct the connected and ordered line sequences (lists of
vertices) $\Sigma_{0i}$
of the adjacent edges of ${\cal F}_0$. The edges are defined by two vertices,
and within $\Sigma_{0i}$ one vertex shares two edges.

{\bf 4.} Construct the connected sequences $\Sigma_i$
from ${\cal F}$ as in {\bf 3}.

{\bf 5.} The projections of the line sequences 
${\cal P}(\omega,\Sigma_{0i})$ cannot intersect each other,
but ${\cal P}(\omega,\Sigma_i)$ can intersect each other and 
${\cal P}(\omega,\Sigma_{0i})$ (intersection of projections can only occur 
when the surface folds away from sight, i.e., the line corresponds to
a facet for which $\langle \nu,\omega \rangle<0$). For any
intersecting projected sequences, find the intersection points $p_{ij}$
on the projection plane with the intersection test above.

{\bf 6.} Define the visible projections $\tilde\Sigma_{0i}$
and $\tilde\Sigma_i$ as the projected sequences 
${\cal P}(\omega,\Sigma_{0i})$ and ${\cal P}(\omega,\Sigma_i)$ that
may have a $p_{ij}$ as an end point.

{\bf 7.} Connect all $\tilde\Sigma_{0i}$ and $\tilde\Sigma_i$ that can form
closed circuits (by systematically comparing
the endpoints of the sequences). The circuit enclosing all the others
(e.g., those due to shadows)
is the approximation of the outer profile contour $\p\cal O$.

If there are more than one generalized profile contours, the 
identification of the circuits should be arranged suitably to enable the
comparison between the model and the data. For example, one shadow region
inside $\p\cal O$ and a smaller $\p{\cal O}_2$ outside $\p\cal O$ due to a
separate closed surface can be identified directly, and all circuits
can be used in determining the best model.

The position of a point in $\p\cal O$ can be parametrized by using
the path length along $\p\cal O$. 
The $\chi^2_\p$ is now given by (assuming one contour per profile)
\be
\chi^2_\p=\sum_{ij} \Vert\varkappa_{\rm o}(c_{ij})
-(\varkappa_{\rm m}(c_{i0}+c_{ij})-\varkappa_{i0})\Vert^2
+\lambda\sum_i (S_i-C_i)^2,\label{aochisq2}
\ee
where o an m stand for observed and modelled,
$0\le c\le 1$ is the normalized path length along the measured
and modelled contours $\p{\cal M}_i$, $\p{\cal O}_i$, 
$\varkappa_{i0}$ is the profile offset for each profile
$i$, $c_{i0}$ is the offset parameter for the path's starting point,
$\lambda$ is a suitable weight factor, and
\be
S_i=\oint_{\p{\cal M}_i} ds,\quad C_i=\oint_{\p{\cal O}_i} ds.
\ee
Thus, $c=s/S_i$ or $c=s/C_i$, where $s$ is the usual path length.

The contour fit can also be modelled by considering the distances of observed
points from the model contour $\p{\cal O}$.
Now we define
\be
\chi^2_\p=\sum_{ij} \inf_s\Big\{\Vert\p{\cal O}_i(s)-
\varkappa_{ij}\Vert^2\Big\},\label{aochisq3}
\ee
where $\varkappa_{ij}$ are the data points,
and we label the points on $\p{\cal O}_i$ by $s$, and assume the
translation due to $\varkappa_{i0}$ to be included in
$\p{\cal O}_i(s)$; here it suffices to consider
points in $\p\cal O$ on whose normal lines $\varkappa_{ij}$ lies. 
When $\p\cal O$ is a set
of line segments, the shortest distance required in (\ref{aochisq3})
(let us denote it by $\delta_{\rm min}$)
is defined by:

{\bf 1.} Let $p'_j\in\R^2$ be the projection of $\varkappa$ on the line
coinciding with the $j$th line segment (corresponding to
$\tan\alpha=(\xi_a-\xi_b)/(\eta_a-\eta_b)$ and 
$(\xi_0,\eta_0)\rightarrow(x,y)$ in the intersection test above):
\be
p'_{(1)}=\frac{D^2\xi-CD\eta+CF}{C^2+D^2},\quad p'_{(2)}=\frac{C^2\eta-CD\xi
+DF}{C^2+D^2}.
\ee 
If the projection is inside the segment,
let $d'_j$ be the distance between $p'_j$ and $\varkappa$.

{\bf 2.} Let $d_k$ be the distance of $\varkappa$ from the $k$th end point 
of the line segments of $\p\cal O$.
 
{\bf 3.} $\delta_{\rm min}$ is the smallest one of all the distances
$d'_j$ and $d_k$. 

In addition to solving for the shape, we usually need to determine the
target's spin state as well in order to have correct projection directions
\cite{kaas01}. In most cases, the target revolves around a constant
pole direction $(\beta,\lambda)\in S^2$ at a constant rotation speed.
The profile plane coordinates $(\xi,\eta)$
are the $x_2'$- and $x_3'$-components of
\be
x'={\sf R}\,x,\label{roteq}
\ee
where
\be
{\sf R}={\sf R}_y(\vartheta-\frac{\pi}{2}){\sf R}_z
(\psi-\lambda){\sf R}_y(-\beta)
{\sf R}_z(-\phi_0-\Omega(t-t_0)),\label{rot}
\ee
where $t$ is the time, $\Omega$ is the rotation speed ($2\pi/P$
for a constant rotation period $P$), $\phi_0$ and the 
epoch $t_0$ are some initial values, and $(\vartheta,\psi)\in S^2$ is
the direction from the target to the observer. We determine
$(\beta,\lambda)$ and $\Omega$ when solving the inverse problem.
It is easy to accommodate 
other spin models such as precession \cite{prec} or nonconstant
rotation speed \cite{yorp} in this formalism.

\subsection{Regularization}

The parameters $P$ describing the target usually have to be 
(moderately) regularized to prevent unrealistic solutions. One aspect
is the smoothness of the body; the larger the target is, the less irregular
it is expected to be. For some parts of the surface this is explicitly 
enforced by the profile data, so the regularization mostly pertains to
the parts covered only by lightcurves that contain little information on
nonconvex features. In those regions, undulation of the surface should be
suppressed (the optimal choice of the suppression weight is discussed in
section 4).

For starlike bodies $\cal B$, a simple (computationally $\chi^2$-like) 
measure of global regularity is
\be
g_S=\int_{\cal B} [r-\langle r \rangle]^2\,d\sigma,
\ee
where $r$ is the radius of the model at the point corresponding to
the surface element $d\sigma$; for polytopes, we can simply use
$g_S=\sum_i (r_i-\langle r \rangle)^2$. 
For such bodies, $g_S$ is typically
quite efficient when the radius is given by a truncated Laplace series
(itself a smoothing agent; see section 5 and \cite{kaast}) and the
regularization weight is low. This is usually the case here
as the profile contours already prevent runaway solutions, so $g_S$
is only needed to polish up the resolution level of shape detail.
For higher weights or models with independent
(uncorrelated) surface points, $g_S$ is not suitable as it
emphasizes global roundedness more than local smoothness.

A measure concentrating on local smoothness (and more suitable
for more complex cases) can be constructed
by considering the negative values of the curvature function.
For polytopes, a practical discrete version of this is
computed by measuring how efficiently the
facets not in the convex hull of the polytope can be blocked (from
viewing or illumination) by their adjacent facets \cite{kaast}. Taking
into account the size and relative tilt angle of the possible blocker
facets adjacent to the facet $i$, 
we can define, e.g., the following measure $\cal C$ 
by summing over the polytope and weighting suitably: 
\be
{\cal C}=\frac{1}{\sum_i A_i}
\sum_{ij} A_{ij}(1-\langle \nu_i,\nu_{ij}\rangle),
\ee
where $A_i$ denotes the area of the facet $i$, and $A_{ij}$ the areas
of those facets around it that are tilted above its plane \cite{kaast} 
(for $i$ in the convex hull, $\sum_j A_{ij}=0$ by definition).
In regularization, we minimize $\cal C$ (for convex bodies 
${\cal C}\equiv 0$).

A further smoothing constraint, to 
be used for non-starlike contours $\p\cal O$ 
if the observations do not cover the profile densely, is given by
augmenting (\ref{aochisq3}) by
\be
\lambda\sum_i\frac{1}{C_i}\oint_{\p{\cal O}_i}
\inf_{\varkappa\in\{\varkappa_{ij}\}}\Big\{\Vert\p{\cal O}_i(s)-
\varkappa\Vert^2\Big\}\, ds,
\ee
which suppresses irregularity on surface parts
not projected near the observed profile points.

A physical constraint for most asteroids is that they are principal-axis
rotators: their maximum moment of inertia is aligned with the rotation
axis due to energy dissipation caused by the nonzero elasticity of the 
material of the body \cite{pravec}. The regularization is defined by
the symmetric inertia tensor \cite{gold}
\be
{\sf I}=\left(\begin{array}{rrr}
P_{22}+P_{33}&-P_{12}&-P_{13}\\
-P_{12}&P_{11}+P_{33}&-P_{23}\\
-P_{13}&-P_{23}&P_{11}+P_{22}
\end{array}\right),
\ee
where the inertia products $P_{ij}$ are
\be
P_{ij}=\int_{\cal B} \rho(x)x_ix_j\,d^3x,
\ee
and here we choose constant density $\rho(x)=1$.
We want to minimize the angle $\tau$ between the $z$-axis of the model and the
eigenvector $I\in\R^3$ (normalized $\langle I,I\rangle=1$) 
corresponding to the largest eigenvalue of
the inertia matrix $\sf I$ of the model shape $\cal B$, 
so we can choose, for example:
\be
g_I=(1-\cos^2\tau)^2=[1-I_3({\cal B})^2]^2,\label{inertreg}
\ee
where the square form $I_3^2$ is useful for weighting purposes and for
removing the sign ambiguity in $I_3$.
A fast way of evaluating $I_3$ in (\ref{inertreg}) for any polyhedron is
described in \cite{dobro}.
Again, profile data constrain the result so strongly that usually the
weight for $g_I$ is very low and sometimes can be set to zero to
obtain, say, $\tau<4^\circ$. Enforcing a $\tau$ much lower than this is 
seldom meaningful due to shape resolution level and inhomogeneities 
in the density.

\section{Optimal combination of data modes: maximum compatibility estimate}

From the statistical viewpoint, 
when we have two or more data modes, we consider their 
simultaneous probability distribution of model parameters and observations
(augmented by prior or regularization 
distributions) in determining the posteriori distribution of the model
and the corresponding estimates. 
The essential problem in this combining is inevitably the 
weighting of distributions. While the
data modes share a common set of parameters describing the underlying
model to be solved for, the models and modalities
of observations may be completely different, and we seldom know a priori
exactly how to compare and weigh their significance.

Let us choose as goodness-of-fit measures (from which probability
distributions can be constructed) the $\chi^2$-functions of $n$ data modes. 
Our task is to
construct a joint $\chi_{\rm tot}$ with well-defined weighting
for each data mode:
\be
\chi_{\rm tot}^2(P,D)=\chi_1^2(P,D_1)+\sum_{i=2}^n \lambda_{i-1}
\chi_i^2(P,D_i)\quad D=\{D_i,i=1,\ldots,n\}
\ee
(to which regularization functions $g(P)$ can be added), where $D_i$ denotes
the data from the source $i$, and 
$P\in\R^p$ is the set of model parameter values.
We assume the $\chi^2_i$-space to be nondegenerate, i.e.,
$$
\arg\min \chi^2_i(P) \ne \arg\min \chi^2_j(P),\quad i\ne j
$$

In two dimensions, denote 
\bea
x(\lambda)&:=&\lbrace\chi_1^2\vert\min\chi_{\rm tot}^2;\lambda
\rbrace,\\\nonumber
y(\lambda)&:=&\lbrace\chi_2^2\vert\min\chi_{\rm tot}^2;\lambda\rbrace. 
\eea
The curve
\be
{\cal S}(\lambda):=[\log x(\lambda),\log y(\lambda)]
\ee
resembles the well-known ``L-curve''
related to, e.g., Tikhonov regularization \cite{belge,engl,hanke}.
However, here we make no assumptions on the shape of $\cal S$.
The curve $\cal S$ is a part of the boundary $\p\cal R$ of the region 
${\cal R}\in\R^2$ formed by the mapping $\chi: \R^p\rightarrow\R^2$
from the parameter space $\pp$ into $\chi_i^2$-space:
$$
\chi=\lbrace\pp\rightarrow
(\log\chi_1^2,\log\chi_2^2)\rbrace,\quad {\cal R}=\chi({\cal P})
$$
where the set ${\cal P}$ includes all the 
possible values of model parameters (assuming that $\chi$ is continuous and
well-behaved such that a connected $\cal R$ and $\p\cal R$ exist). 
If the possible values of $\chi^2_i$ are not bounded, the remaining part 
$\p\cal R\setminus\cal S$ stretches droplet-like towards $(\infty,\infty)$.
The parameter $\lambda$ describes the
position on the interesting part ${\cal S}\subset\p\cal R$,
and it is up to us to define a criterion
for choosing the optimal value of $\lambda$.
 
The logarithm
ensures that the shape of ${\cal S}(\lambda)$ is invariant under unit or scale
transforms in the $\chi_i^2$ as they merely translate ${\cal S}$ in the
$(\log\chi_1^2,\log\chi_2^2)$-plane. It also provides a meaningful metric
for the $\log\chi_i^2$-space: distances depict the relative difference
in $\chi^2$-sense, removing the problem of comparing the absolute values of 
quite different types of $\chi_i^2$.
The endpoints of ${\cal S}(\lambda)$ are at
$\lambda=0$ and $\lambda=\infty$, i.e., at the values of $\chi_i^2$
that result from using only one of the data modes in inversion. We can
translate the origin of the 
$(\log\chi_1^2,\log\chi_2^2)$-plane to a more natural
position by choosing the new coordinate axes to pass through these endpoints.
Denote
\bea
\hat x_0&=&\log x(\lambda)\vert_{\lambda=0}=\log\min\chi_1^2\\\nonumber
\hat y_0&=&\log y(\lambda)\vert_{\lambda\rightarrow\infty}=\log\min\chi_2^2. 
\eea
Then the ``ideal point'' $(\hat x_0,\hat y_0)$ is the new origin in the
$(\log x,\log y)$-plane. A natural choice for an optimal location
on $\cal S$ is the point closest to $(\hat x_0,\hat y_0)$, i.e., the
parameter values $P_0\in\pp$ such that
\be
P_0=\arg\min\Big([\log\chi^2_1(P)-\hat x_0]^2+[\log\chi^2_2(P)-
\hat y_0]^2\Big),
\ee
so we have, with $\lambda$ as argument,
\be
\lambda_0=\arg\min \Big([\log x(\lambda)-\hat x_0]^2+[\log y(\lambda)-
\hat y_0]^2\Big).\label{lambda0}
\ee
In this approach, neither the numbers of data points
in each $\chi^2_i$ nor the noise levels as such affect the
solution for the optimal $P_0$ as their scaling effects cancel out in
each quadratic term. $P_0$ is thus a pure compatibility estimate describing
the best model compromise explaining the datasets of different modes
simultaneously.

We call the point $P_0$ the
{\em maximum compatibility estimate} (MCE), and $\lambda_0$ the {\em maximum
compatibility weight} (MCW). This corresponds to the maximum
likelihood estimate in the case of one data mode, or to the maximum a
posteriori estimate as well since we can
include regularization functions here. If regularizing is used, the weights
for the functions are either determined in a similar manner (see below),
or they can be fixed and the regularization terms are absorbed in
$\chi_1^2$ (otherwise ${\cal S}\subset\p\cal R$ does not hold).

Another choice, frequently used in the L-curve approach, is to find the
$\lambda$ at which $\cal S$ attains its maximum curvature \cite{hanke,engl}, 
but evaluating this point is less robust than finding $\lambda_0$,
and (\ref{lambda0}) is a more natural prescription, requiring no assumptions
on the shape of $\cal S$. We make two implicit assumptions here:
\begin{enumerate}
\item
The solutions $P_{\p\cal R}$ 
corresponding to points on $\p\cal R$ should be continuous (and one-to-one) 
in $\pp$-space along $\p\cal R$ at least in the vicinity of the solution
corresponding to $\lambda_0$. If this is not true (in practice, 
if $P_\lambda=\arg\min\chi^2_{\rm tot}(P)$
makes large jumps in $\pp$ for various $\lambda$ around $\lambda_0$), 
one should be cautious
about the uniqueness and stability of the chosen solution $P_0$, and
restrict the regions of $\pp$ included in the analysis. 
\item
The optimal point $\lambda_0$ on 
$\cal S$ should 
be feasible: if we have upper limits $\epsilon_i$ to acceptable $\chi^2_i$,
the feasible region $\cal F$
is the rectangle $\bigcap_i\lbrace\log\chi^2_i\le\log\epsilon_i\rbrace$. If
$[\log\chi_1^2(P_0),\log\chi_2^2(P_0)]\notin
\cal F$ and ${\cal F}\cap{\cal R}\ne
\emptyset$, we choose the point on the portion ${\cal S}\subset\cal R$ 
closest to the one corresponding
to $\lambda_0$ (i.e., $\log\chi^2_i=\log\epsilon_i$ for one $i$). 
If ${\cal F}\cap{\cal R}=\emptyset$, the data modes do not
allow a compatible joint model, so either the model is incorrect for 
one or both data modes, or one or both $\epsilon_i$ have been estimated 
too low (e.g., systematic errors have not been taken into account). Note
that model insufficiency should be taken into account in the estimation of 
$\epsilon_i$.
\end{enumerate}

Note that, in the interpretation ${\cal R}=\chi({\cal P})$,
$\lambda$, $\chi_{\rm tot}^2$ and $\p\cal R$ are all
in fact superfluous quantities, and we can locate the point estimate
MCE $P_0$ entirely without them with standard optimization procedures
(and with no extra computational cost).
However, it is useful (though computationally
somewhat noisier) to approximate $\cal S$ via the minimization 
of $\chi_{\rm tot}^2$ with sample values of $\lambda$ (see Fig.\ 1), 
as in addition to obtaining the MCW $\lambda_0$ (and hence MCE as well) we can
plot $\cal S$ to examine the mutual behaviour of the complementary data
sources (including the position of the feasibility region $\cal F$ 
w.r.t. $\cal S$). The solution for $\lambda_0$ is also needed for
constructing distributions based on $\chi^2_{\rm tot}$. Another possibility
to examine $\cal R$ and $\p\cal R$ is direct adaptive Monte Carlo sampling, 
but this is computationally slow.

This approach straightforwardly generalizes to $n$ $\chi^2$-functions and
$n-1$ parameters $\lambda_i$ describing the position on the $n-1$-dimensional
boundary surface $\p\cal R$ of an $n$-dimensional domain $\cal R$: the MCE is
\be
P_0=\arg\min\sum_{i=1}^n\Big[\log\frac{\chi^2_i(P)}{\chi^2_{i0}}\Big]^2,
\quad\chi^2_{i0}:=\min\chi^2_i(P),\label{p0}
\ee
and the MCW is
\be
\lambda\in\R^{n-1}:\quad \lambda_0=\arg\min\sum_{i=1}^n\Big[\log
\frac{\hat\chi_{i,{\rm tot}}^2(\lambda)}{\chi^2_{i0}}\Big]^2,
\quad\hat\chi_{i,{\rm tot}}^2(\lambda):=
\Big\{\chi^2_i\Big\vert\min\chi^2_{\rm tot};\lambda\Big\}.\label{l0}
\ee

Another scale invariant version of MCE can be constructed by plotting 
$\chi^2_i$ in units of $\chi^2_i/\chi^2_{i0}$ and shifting the new origin
to $\chi^2_i/\chi^2_{i0}=1$:
\be
P_0=\arg\min\sum_{i=1}^n\Big[\frac{\chi^2_i(P)}{\chi^2_{i0}}-1\Big]^2,
\quad\lambda_0=\arg\min\sum_{i=1}^n\Big[
\frac{\hat\chi_{i,{\rm tot}}^2(\lambda)}{\chi^2_{i0}}-1\Big]^2.
\ee
This, however, is exactly the first-order approximation of (\ref{p0}) and
(\ref{l0}) in $\delta\ll 1$ when $\chi^2_i/\chi^2_{i0}=1+\delta$, giving
virtually the same result as (\ref{p0}) and (\ref{l0}) as usually
$\chi^2_i(P_0)/\chi^2_{i0}-1\ll 1$ in the region around $\chi^2_i(P_0)$, 
and any larger ratios of $\chi^2_i/\chi^2_{i0}$ are not eligible for
the optimal solution (see Fig.\ 1).

Instead of the $L_2$-norm $\chi^2$ (and
the corresponding $\chi^2$-distribution), we can
choose some other goodness-of-fit measure $\varepsilon(P,D)\ge 0$
(and distribution) for the individual data modes. For a linear combination 
of these, we have
$$
\varepsilon_{\rm tot}(P,D)=\varepsilon_1(P,D_1)+\sum_{i=2}^n 
\lambda_{i-1}\varepsilon_i(P,D_i).
$$
In lightcurve measurements, for example, the effect of systematic 
errors in both model and data dominates over random noise when
the noise level is not high
\cite{iau}, so it is not mandatory to use $\chi^2$ as a standard
measure of fit.

\begin{figure}
\begin{center}
\includegraphics[width=11cm]{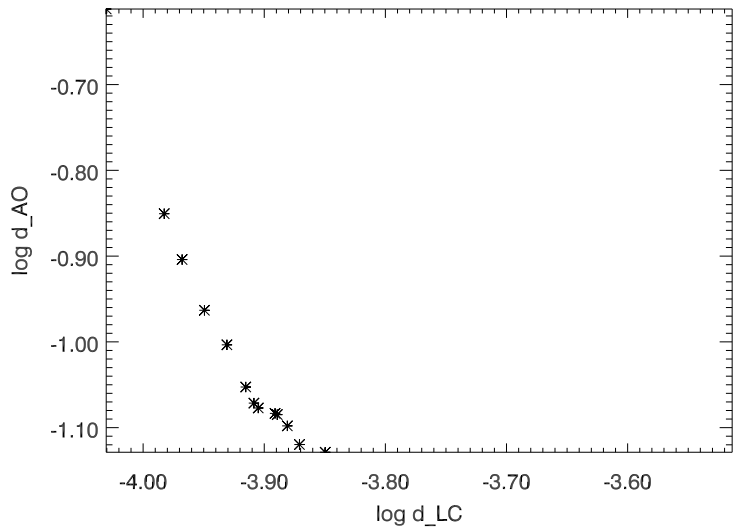}
\caption{$\cal S$ curve plotted for 2 Pallas with various weights $\lambda$
(LC for lightcurves, AO for adaptive optics profiles).}
\label{scurve}
\end{center}
\end{figure}

It is possible to use this approach for general regularizing functions 
$g(P)$ as well (change $\chi_i^2\rightarrow g(P)$ for some $i$), 
but in such cases the shape of $\cal S$ must be taken into account.
If it is possible to have a solution $g(P')=0$ for a regularizing function
$g$ (or an almost vanishing $g(P')$ such that 
$\log g(P')\rightarrow -\infty$), 
the above scheme automatically returns $P'$ and ignores the actual 
data altogether. Thus one should, e.g., set a lower practical limit to 
$g(P)$ by looking at the shape of $S$, and choose the $\lambda_0$ within
the restricted part of $S$. Likewise, one can use the above scheme for
assigning noise-level-independent 
weights to subsets of the same data mode (rather than have
the standard $\chi^2$ evaluated from all data points), but obviously the
subsets cannot be chosen arbitrarily if the result is to make sense. For 
example, one can estimate the optimal weight for one lightcurve that
appears to reveal features not contained in other lightcurves and thus
judge its real significance. Even one noisy lightcurve with a few points, 
taken at a special observing geometry, may well contain significant 
information that needs to be weighed more against less noisy but more 
ordinary lightcurves.

\section{Numerical implementation}

As examples of the optimal combining of lightcurves and AO profiles, we show
representative results for the asteroids 2 Pallas and 41 Daphne. 
Full detailed descriptions of 
the observations and models of these targets
are presented in \cite{carry2} and Carry et al. 
(in preparation). An example of an even more irregular shape
constructed with our procedure is the
model of the primary body of the binary asteroid 
121 Hermione \cite{descamps}.
The lightcurve $\chi^2_L$ was computed as in \cite{kaast,kaas01} and
profile $\chi^2_\p$ as in the starlike case of
section 3, and the minimization of 
$\chi^2_{\rm tot}$ was performed as in \cite{kaast,kaas01}.
The observed profiles are projections of the target on the plane-of-sky
$S^2$ converted to pixels on the instrument plane, while the model is
constructed in absolute (km) size, so the model/profile scale conversion
is given by the AO instrument's angular resolution and the distance
between the target and the observer.
The profile contour extraction procedure with 
wavelets (as an average of several AO images obtained
 in a short time interval) is described in 
\cite{carry1,carry2}. 

In general, the resolution
of the model must be somewhat lower than the apparent resolution of the
AO images as the sparse profile samples will produce artificial features
elsewhere in the model if a near-perfect profile fit is enforced (even if the
observed profile details were exactly right). The inverse problem
has thus some ill-posedness at local scales starting near the profile 
resolution level, but the ill-posedness at more global scales, inherent
to lightcurve data \cite{kaas01,genproj,iau}, is removed with AO profiles.
The weight factor $\lambda$ mostly
takes care of this, and fine-tuning is obtained with $\lambda_S$ for
the smoothness constraint $g_S$. 
For the examples here, the weight of the inertia regularization
function $g_I$ was low as there were several profiles available;
virtually the same result was achieved with $\lambda_I=0$.
The weights $\lambda$ and $\lambda_S$ were 
determined with the scheme of section 4; the examined interval of
$\lambda_S$ was restricted to realistic values corresponding to the
resolution level of the AO images. 

Fig.\ 1 depicts a typical evaluation
of the curve $\cal S$ at various choices of $\lambda$; or rather, this plot
portrays the cross-section of the 2-surface $\p\cal R$ in $\R^3$ with
$\lambda_S$ fixed at its final optimal value.
The values for $\chi_i^2$ are normalized to be the rms deviations
of model fits $d_i=\sqrt{\chi_i^2/N_i}$, as in logarithmic scale
this corresponds only to a shift of origin and a uniform linear change
of plot scaling. The plotted points outline 
the curve ${\cal S}(\lambda)$ 
that is rather an oblique line than an L-shape, and the ideal
point region, i.e., the point closest to the lower left-hand corner,
can directly be found. The endpoints $\lambda=0$ and
$\lambda=\infty$ stop at saturation regions rather than continue to
large distances in the $\log\chi^2$-space.
As can be seen from Fig.\ 1, 
computational noise in the estimated points at
$\lambda=0$ and $\lambda=\infty$, corresponding to a small change of the 
position of the new origin w.r.t. $\cal S$, 
does not affect the estimated location
of the optimal point on $\cal S$ significantly.

\begin{figure}
\begin{center}
\includegraphics[width=5.8cm]{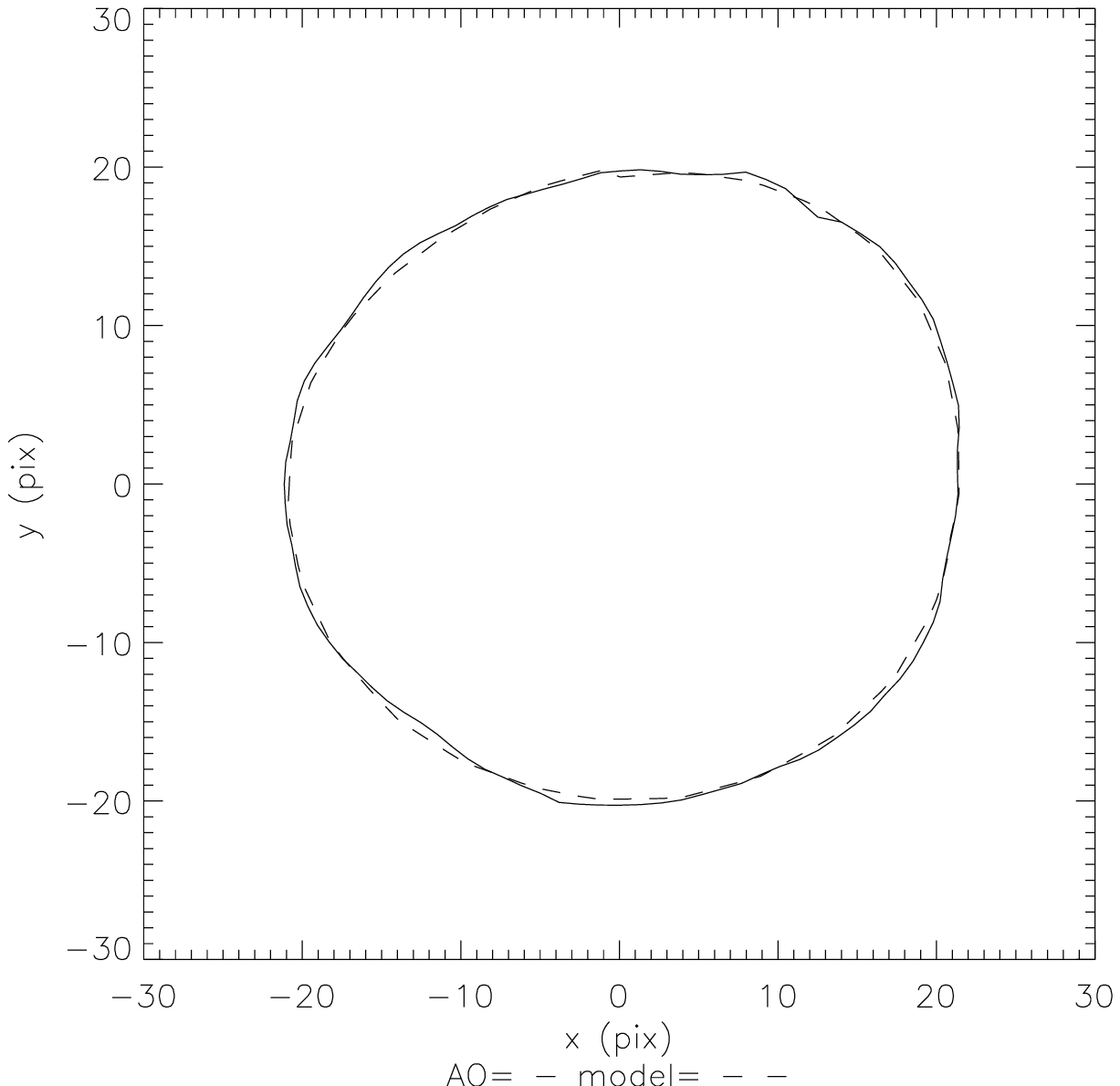}\includegraphics[width=5.8cm]{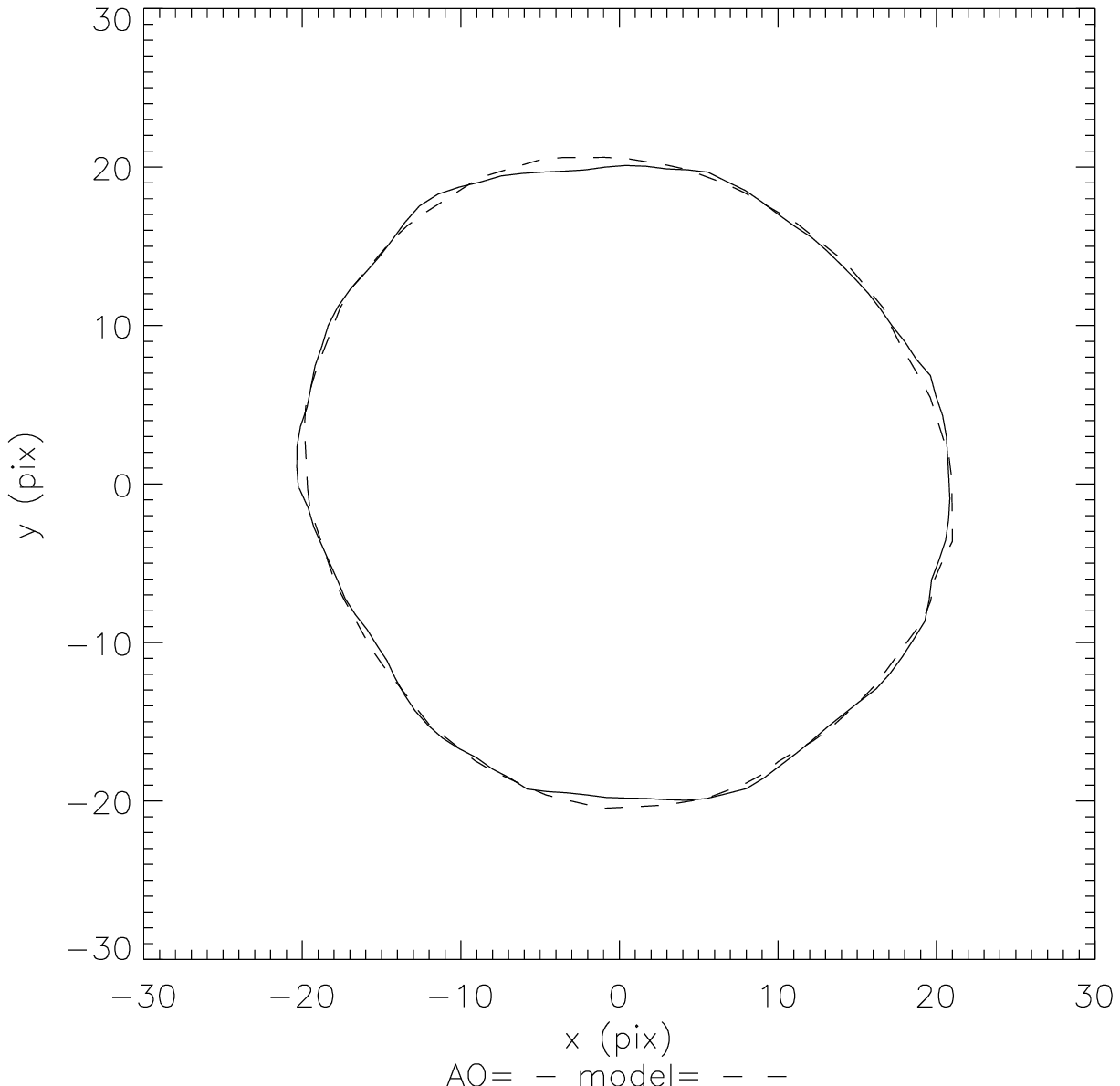}\includegraphics[width=5.8cm]{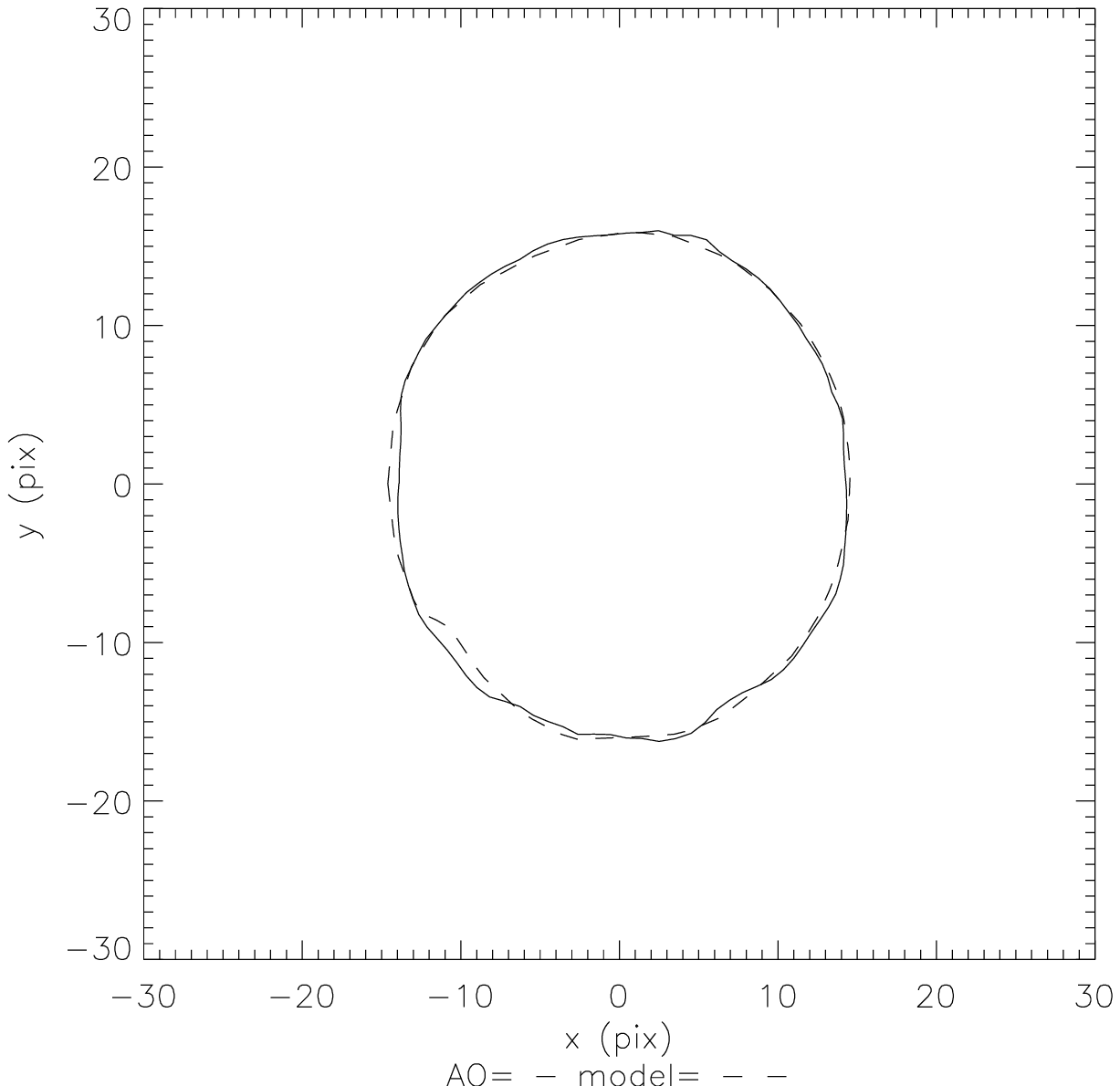}
\caption{Sample observed (solid lines) vs. modelled (dashed lines)
AO contours for 2 Pallas. Coordinates are in pixel units.}
\label{aopallas}
\end{center}
\end{figure}

Sample observed vs.\ modelled profiles for 2 Pallas and 41 Daphne are shown
in Figs.\ 2 and 3.
The starlike surface model was described by the exponential
Laplace (spherical harmonics) series for the surface radius $r$ \cite{kaast} 
\be
r(\theta,\varphi)=\exp\Big[\sum_{lm} c_{lm}Y_l^m(\theta,\varphi)\Big],
\quad (\theta,\varphi)\in S^2,
\ee
truncated at suitable $l,m$, with $c_{lm}$ as the shape parameters 
to be solved for. Other model parameters are the profile offset 
$(\xi_0,\eta_0)$ for each image and the spin parameters.
For asteroid 2 Pallas (a rather spherical body with size class 500 km), 
the Laplace series was truncated
at maximal $l=6,m=6$, while for the more irregular
41 Daphne (size class 200 km)
the truncation point $l=8,m=6$ was more appropriate.
The early truncated Laplace series and the choice
of the truncation point are implicit regularization measures as such. We
leave the discussion of the choice of model discretization level elsewhere 
(cf.\ \cite{kaipio}) as here its effect on the data mode weighting 
was negligible (within a feasible set of choices), and the resolution level
of AO images (as well as keeping $\lambda_S$ low and avoiding
artificial surface features) essentially
determined the choice in practice after some sampling.
For AO data, the choice of the Laplace series as a model 
is practical, while for, e.g., 
detailed space probe data a mesh of independent surface points is more
accurate and computationally feasible.

\begin{figure}
\begin{center}
\includegraphics[width=5.8cm]{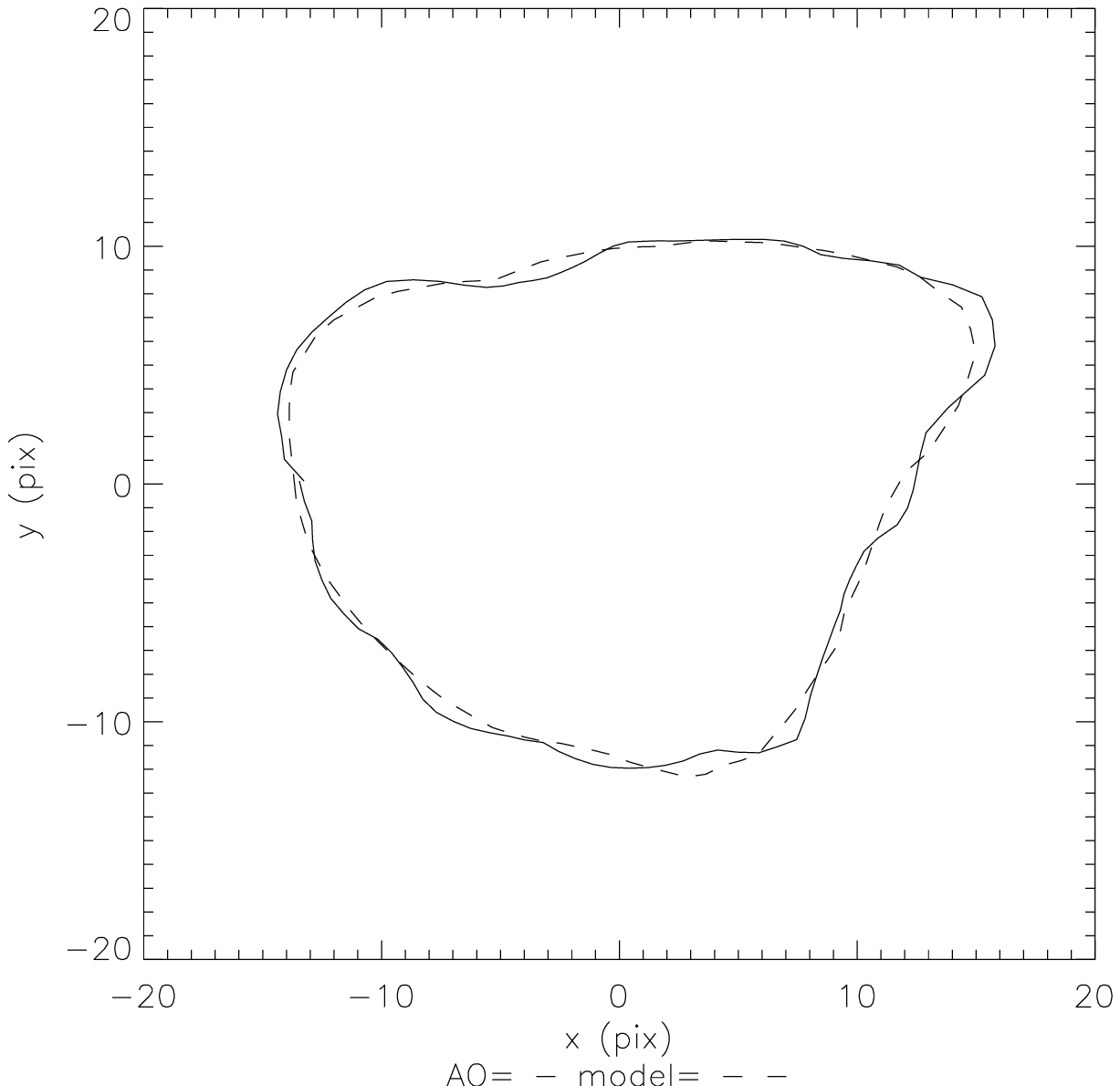}\includegraphics[width=5.8cm]{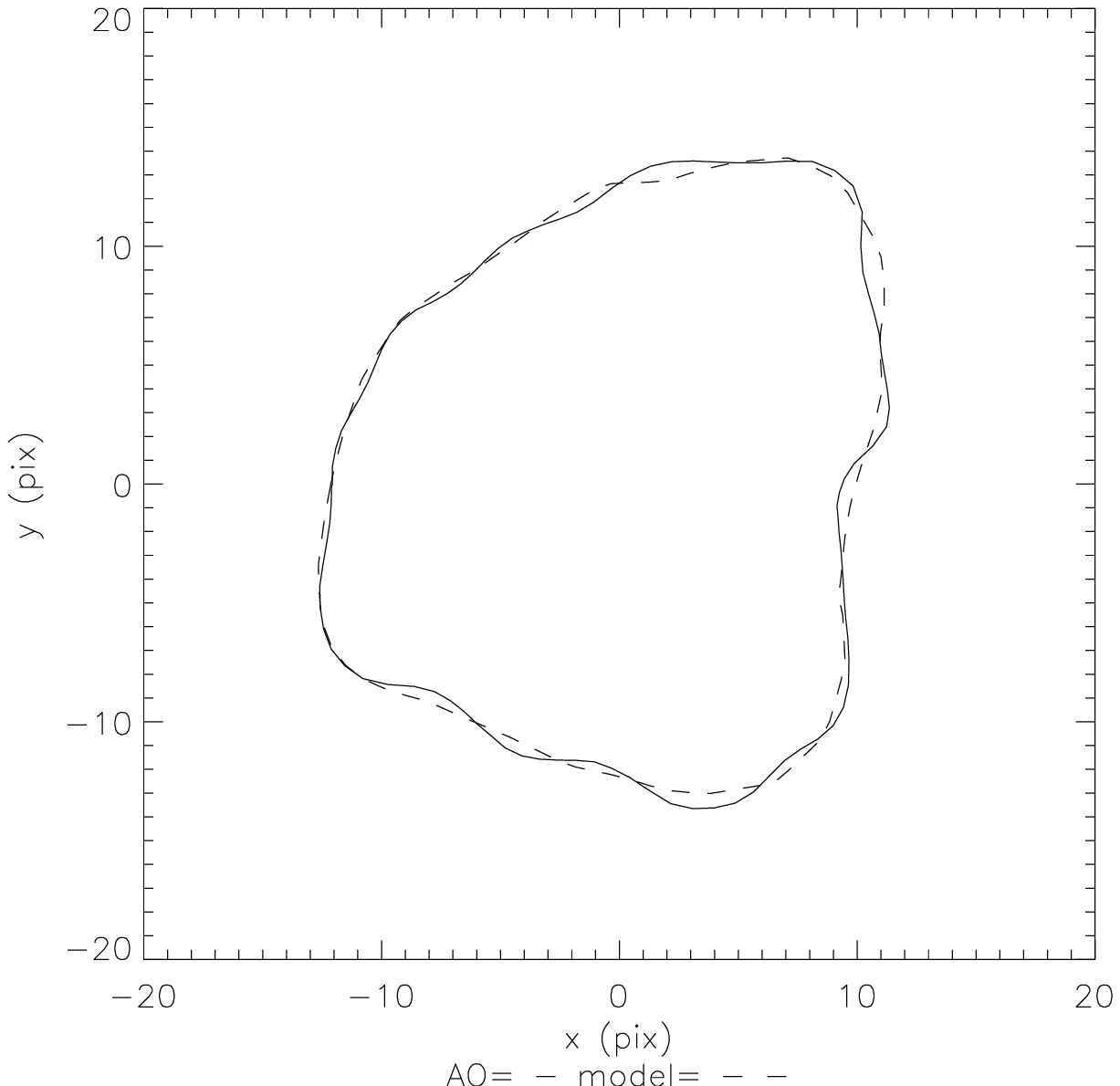}\includegraphics[width=5.8cm]{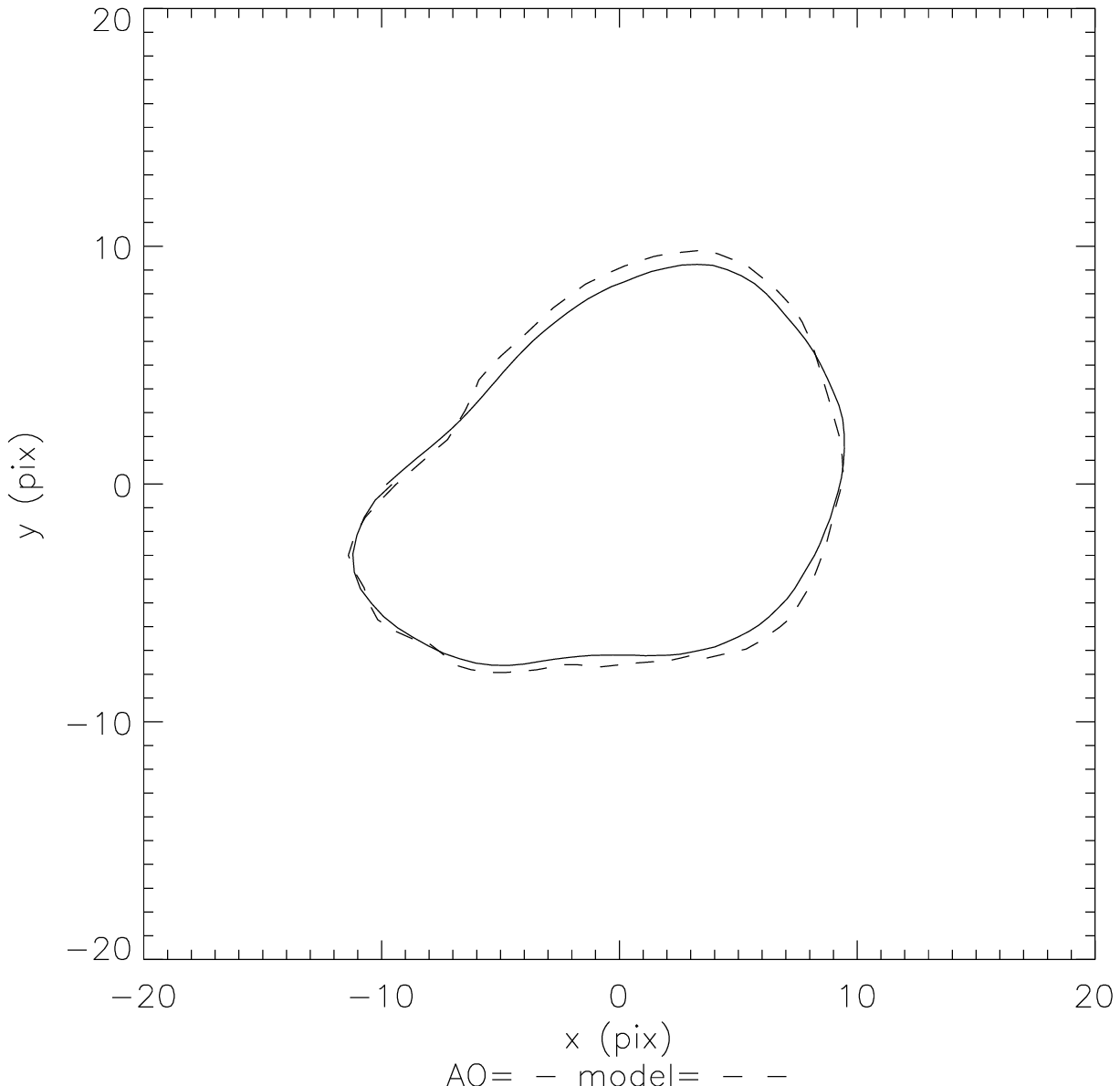}
\caption{Sample observed (solid lines) vs. modelled (dashed lines)
AO contours for 41 Daphne. Coordinates are in pixel units.}
\label{aodaphne}
\end{center}
\end{figure}

Once the weight factors $\lambda$ and $\lambda_S$
are determined, the result is usually 
stable and restricted to one region in the parameter space $\pp$:
probing feasible solutions $P$ corresponding to $\chi^2_{\rm tot}(P)$ slightly
lower than $\chi^2_{\rm tot}(P_0)$ produces essentially the same results.
Due to restricted orbital geometries, lightcurve data alone often 
imply two almost equally possible pole directions with mirror-like shape
solutions \cite{genproj,iau}; even one AO (or other) image 
usually resolves this typical ambiguity \cite{keck}. The result is also
typically stable w.r.t. weights in the vicinity of MCW. The obtained
MCE appears to be well justified when one samples the solutions along
$\cal S$: it provides a very good match to profile details without
straying far from the observed lightcurves, and does not predict too
prominent features on the parts of the surface not projected onto the
profile contours.

\section{Conclusions and discussion}

We have examined the classes of shapes reconstructable 
by the (generalized) profiles of objects in $\R^3$,
and presented a method for using lightcurves and
the observed contours of generalized profiles simultaneously to
produce shape (and spin) models  
with more details (and a lower degree of ill-posedness)
than in the pure lightcurve mode. We have also shown that there is
a well-justified criterion and an efficient method for determining the
optimal weighting of data modes. Applied to real data, the method works
very well, and we can use simple regularization functions. 
In addition to adaptive optics observations, asteroid
profiles can also be obtained from other sources such as interferometry,
space telescopes, and stellar occultations (partial profiles).

The use of profiles is practical as it removes two sources of systematic
errors inherent to using full images (brightness distributions $\cal I$ on the
image plane): the errors in $\cal I$ from AO deconvolution and the model
$\cal I$ errors due to the insufficently modellable light-scattering 
properties of the surface of the target body. On the other hand,
profile determination requires the data to be sharp enough, not with fuzzy
images. If the images are fuzzy, we usually have to resort to using some
brightness and blurring model for fitting full images, even though the result
will be less certain.

The concept of the maximum compatibility estimate is directly applicable
to any inverse problems with complementary data modes. The invariance 
properties of the MCE make it more generally usable than heuristic strategies
for choosing the weights, especially when they use assumptions on the
shape of $\p\cal R$ or other case-specific characteristics.

\subsection*{Acknowledgements}

It is a pleasure to thank Benoit Carry and Josef \v{D}urech 
for discussions and comments. The sample adaptive optics data used
in figures here are courtesy of B. Carry, A. Conrad, J. Drummond, 
C. Dumas, S. Erard, and W. Merline. This work was supported by the Academy
of Finland (project ``New mathematical methods in planetary and
galactic research'').


\medskip
{\it E-mail address:} First.Lastname [at] tut.fi

\end{document}